\documentclass{amsproc}

\theoremstyle{definition}

\theoremstyle{remark}

\numberwithin{equation}{section}

\begin{document}

\title{Iron Rings, Doctor Honoris Causa Raoul Bott,
\\ Carl Herz, and a Hidden Hand} 

\author{P. Robert Kotiuga}
\address{Department of Electrical and Computer Engineering, Boston University, Boston, MA 02215}
\curraddr{Department of Electrical and Computer Engineering,
Boston University, Boston, MA 02215}
\email{prk@bu.edu}

\date{June 18, 2009.}


\begin{abstract}
The degree of Doctor of Sciences, honoris causa, was conferred on
Raoul Bott by McGill University in 1987. Much of the work to make this
happen was done by Carl Herz. Some of the author's personal
recollections of both professors are included, along with some context
for the awarding of this degree and ample historical tangents. Some cultural
aspects occurring in the addresses are elaborated on, primarily, the
Canadian engineer's iron ring. This paper also reprints both the
convocation address of Raoul Bott and the presentation of
Carl Herz on that occasion.\thanks{An edited and reformatted
  version of this paper, with an additional photo, will appear in
  a volume dedicated to Raoul Bott\cite{K10}. The author hopes to
  expand on some aspects of this preprint in future versions.}
\end{abstract}

\maketitle

\section*{Introduction }

Raoul Bott needs no introduction in this volume. However, reprinting
 his address at the 1987 McGill convocation both gives some insight 
into the effort to award him an Honorary Doctorate in Mathematics from
 McGill, and a context to develop some less than mathematical themes,
 unashamedly from the point of view of an electrical engineer who
 enjoys the historical aspects of his discipline. Early 
on I was tipped off that Carl Herz was behind the effort, and the
 memory of Prof. Herz made me realize that I had to
 follow the trail like a hound. The result complements more
 technical presentations and anecdotes pertaining to
 Montreal\footnote{For anecdotes from Bott's years in Montreal see
 \cite{Tu} and Candace Bott's remarks in this volume.}. 
I am grateful to Candace Bott for digging up her father's
 commencement\footnote{Although ``convocation'' and
 ``commencement'' have different meanings in general, in the context of
 a university graduation ceremony they have the same meaning in Canada
 and in the USA respectively. For the purposes of this article, they
 are used interchangably.} address, and to Dominique 
 Papineau of McGill University who showed up in Boston with a
 complete file pertaining to the awarding of the honorary degree from
 McGill's archives. The support of various McGill faculty who 
 listened to me think aloud is much appreciated, namely Peter Caines,
 Jacques Hurtubise, Joachim Lambek, and Peter Russell.

\section*{How I got to know Raoul Bott}

I can't remember who first connected Raoul Bott to McGill Engineering
 in my mind; most likely it was Peter Caines, Robert Hermann or Carl
 Herz. However, I do remember
 attending his talk in the Physics department at McGill in, I
 believe, 1982. These days, if the fact that I dragged my wife to
 be to the talk is mentioned, and that she  
 coped with my enthusiasm with a sense of humor, my kids will kindly remind
 me of all my thrifty ideas for a good time! During those years I spent
 many hours ``teething'' on the book of Bott and Tu, and as a
 NSERC post-doc in the MIT Mathematics department in 1985 I finally
 got the opportunity to audit a course given by Raoul Bott.  Not only
 was his mathematics entrancing, but we were both McGill
 engineers!  During the first lecture he spotted my iron ring. I
 introduced myself afterwards, and we talked about a variety of things.  Before
 long, I made a habit of auditing every course he taught. About a year
 later, he told me that he would be the recipient of an honorary
 doctorate from McGill --- forty years after McGill wouldn't have him as a
 graduate student in the mathematics department (he was unwilling to complete a
 second undergraduate degree in mathematics). Past history aside, he
 seemed genuinely honored, but he didn't quite know who to share this
 news with. Reading the convocation speech, I now see it as his way to
 make peace with history and a means to repackage it constructively
 for graduates forty years younger than himself.  Stepping back from the
 ceremony from decades ago, the reader is invited to read his
 ``Autobiographical Sketch''\cite{BoA}.

I clearly enjoyed his lecturing style as well as the attention and
 discipline he demanded in the classroom. It certainly
 contrasted with the convocation speech. In the first lecture, he
 encouraged students to ask questions and claimed that he liked
 ``stupid questions'' because he could answer them. He was very serious
 about this and if we didn't realize it initially, we eventually
 learned that he put the bar very high, and even higher for himself:
 His answers to questions were always more profound than the original
 questions and he once walked out of his own lecture in frustration
 because he didn't like it! That was dramatic, baffling, and
 unexpected --- especially since we knew he liked to think on his
 feet. In this particular case, he reappeared the next class, his
 arguments were impeccably clear and elegant, and his credibility was
 not only 
 restored, it soared! Going beyond the classroom, Bott's mentoring
 of graduate students is the source of legend. This volume has plenty
 of testimony from his students, and the reminiscences by Robert
 MacPherson \cite{MP} testify to his high expectations. The convocation
 speech is very different in that it demonstrates his ability to
 connect with an audience that has never seen him and most likely will
 never see him again. Before attempting to advise graduates, he warms
 up the audience by telling of the pivotal event in his professional
 life and saying: ``I tell you all this only in part as a jest.'' Only when
 the stage is set does he give the essence of his address credibly and
 in a few sentences. The message 
 appears in a flash after saying ``But my time is up!''. The style
 mirrors his approach to giving a colloquium talk. 
 
Bott enjoyed cultivating certain habits which were best left
unmentioned in the convocation speech. For instance, he lectured at 8:30
a.m. in order to have a flexible day after 10:00 a.m.! Graduate
students felt this cramped both their style and sleeping habits, and
this is where some of Bott's Old World sensibilities kicked in when they
dared to doze off in class. He could toss chalk and have it land on
the table inches from the sleeping student's face, startling them. He
clearly relished doing so.  The memorable line which would meet the
startled face changed with every successive offense. From ``didn't want
you to miss anything'' to ``how much is tuition at Harvard?'' to ``who pays
your tuition?'' In Old World style, this was all for the benefit of the
student and there was little room for self-preservation. These days it
isn't easy for a professor to be respected for doing this in a private
university where students can feel like paying customers once the
tuition bill is settled. Somehow Bott was 
consistently more mischievous than the students, and got away with
it. Clearly he had extensive experience testing his teachers
and this experience always gave him the upper hand in the classroom.

As for the impact Bott's lectures made on me, I'd be treading on thin
ice if I tried to say why they were fantastic.  Loring Tu says
that\footnote{See ``Reminiscences of
Working with Raoul Bott" in \cite{STY}} Victor Guillemin, at a
conference celebrating his 60th birthday, proudly announced that he
took twelve courses from Bott, and to Loring's chagrin, he could only
list eleven. Clearly, I am in no position to speak with authority
about Bott's lectures! In my case, I loved his lecture style, the
lecture material, and I felt a definite kinship since I could ask
tangentially related questions after class and consistently get
profound answers. There is perhaps one personal anecdote I can add to
the many that I've heard.  One day in class, after Bott explained the
set-up of the Lefschetz fixed point theorem in terms of the transverse
intersection of the graph of a map from a manifold to itself with the
graph of the identity map, he claimed that, by duality, the Lefschetz
number could be easily computed by picking a basis for integral
cohomology, pulling back by the appropriate maps, taking wedge
products with Poincare duals and integrating. When he claimed that it
reduced to basic matrix algebra involving the induced automorphisms on
cohomology groups which an engineer could do, the class just didn't
make the type of eye contact he was hoping for. At that point he
called me up to the board and told me to fill in the details of the
calculation! As I (methodically) wrapped up the calculation, he
identified me as an engineer, emphasized that budding topologists
shouldn't shy away from such concrete calculations, and took
satisfaction in the fact that he made his point. In my mind he
reinforced the fact that Daniel Quillen's thesis advisor could
make us think functorially while, as a student of Richard Duffin, he
could encourage us to ``think with our fingers'' and to always
maintain a balance between the conceptual and the computational.

Obviously, I'm enamored with Raoul Bott, and thrilled that he was the
recipient of an honorary degree from McGill. However, my purpose here
must be more focused. Specifically, in my mind, a few key points need
elaboration: 
\begin{itemize}
\item 	Who was the driving force behind getting McGill to award Bott
an honorary degree in Mathematics forty years after he left McGill as
an engineer? Clearly, Bott was deserving, but it takes a kindred
spirit to overcome the inertia of a bureaucracy and Carl Herz was such
a kindred spirit.
\item In his convocation address, Bott vividly describes the pivotal
  moment at McGill when he decided to become a mathematician. 
  However, how he was going to do it was not at all clear at
  the time. The wonderful and profound connections between topology
  and physics have been studied intensely in recent decades, but what is
  needed is a hint of the path from ``engineering mathematics'' to the
  mathematics Bott is known for\cite{Bo85},\cite{JB}. In retrospect,
  it almost seems that 
  this path could have been more clear to Gauss and Maxwell than to
  modern specialists. We'll soon see, the career of Hermann Weyl
  provides us with a perspective and some key insights.
\item I'm fascinated with Bott's struggle to reconcile old and new
world sensibilities.  I saw this encoded in my
interactions with him and in the convocation
address. I use the word ``encoded'' in reference to the address, because
his references to the uniquely Canadian iron ring have their roots in
various Quebec City bridge disasters, Rudyard Kipling's poem ``Sons of
Martha,'' and associated Biblical references. These details are
required to fully decode the message. 
\end{itemize}

\section*{Carl Herz as a gateway to history}

I can distinctly remember the day Carl Herz knocked me off my feet.
At the time I was a graduate student in Electrical Engineering at
McGill and he was a feisty and famous professor of Mathematics. We
began to chat after some seminar in the EE department, and he asked me
what I was doing for a thesis. He listened as I told him how I felt that
the reformulation of Maxwell's equations in terms of differential
forms was essential for the resolution of some key problems in
computational electromagnetics. Specifically, most of the boundary
value problems in low frequency electromagnetics amounted to Hodge
theory on manifolds with boundary, with the periods of harmonic forms
identified with the variables found in Kirhhoff's laws. Furthermore, I
told him that the whole framework has a variational setting which can
be discretized by appealing to ``Whitney forms'' in order to obtain a
finite element discretization with desirable properties. To me it was
all obvious if one read the papers written by Donald Spencer and his
students in the 1950's and interfaced them with Whitney's ``Geometric
Integration Theory''. In retrospect, this was a natural connection
given the work of Jozef Dodziuk and Werner Muller's proof of the
Ray-Singer conjecture, but it was not apparent at the time. Carl
listened, started pacing back and forth and 
I was beginning to worry that he was going into a trance! I don't know
what was going on in his mind, but I braced myself for what could come
out of his mouth.

I think I was standing in stunned silence when Carl stopped and
asked me if I ever read Maxwell. Sheepishly, I told him that I read a
good deal of Maxwell and that everyone in my field swears by
Maxwell. He then asked me if I knew what a periphractic number
was. When I expressed my ignorance, he went on to point out that
Betti's paper was written a decade after Maxwell's treatise, and that
in Maxwell's treatise the first Betti number was called a ``cyclomatic
number'' - a term introduced by Kirchhoff, and still used in graph
theory. He went on to tell me that the second Betti
number was called a ``periphractic number''... I later found out that
Maxwell borrowed the term from
Listing\footnote{See Breitenberger \cite{Br} in James \cite{J} for an
article on Johann Benedikt Listing 
and his book.} and that Listing was the person who coined the term
topology. In one swoop Carl convinced me that Maxwell was often quoted
but never read, and that if I wanted to get to the origin of these
topological ideas the origin would be in some language other than
English. Clearly, I was humbled --- but I felt better when I looked up
``periphractic'' in the unabridged Oxford dictionary and found that
Maxwell's treatise is the first and last use of the word in the
English language!\footnote{It is irresistable to point out the
  connection between Maxwell and Morse theory in this article about
  Bott. Listing\cite{L} is credited as being the first to
  systematically obtain a cell decompositions of 3-manifolds by
  tracking the change in topology as level sets cross a 
critical point. Maxwell \cite {M70} then wrote a paper citing Cayley and
Listing. In his treatise\cite{M}, Maxwell uses the rudiments of Morse
theory with the fact that a harmonic function cannot achieve a maximum or
minimum in the interior of a region in order to make topological deductions.}

Besides being awed by Carl's encyclopedic knowledge, there are two big
lessons I have leaned over the years and which were initiated by my
encounter with Carl and other mathematicians from his generation who
``read the masters.'' The first was that Maxwell had a profound
experimental and 
theoretical knowledge, and that much of the inspiration for his
theoretical work came from reading and corresponding with Germans
(Gauss, Riemann, Kirchhoff, Clausius, Helmholtz,
Listing, ...). Furthermore, it was the Germans who took Maxwell
seriously when no one else did- from Helmholtz' student Hertz
demonstrating radio waves, to Boltzman developing statistical
mechanics, and to Einstein developing the logical physical conclusions of
Maxwell's theory. Contrast this with the situation in England where
Oliver Heaviside was  
considered a self-educated eccentric who died in poverty despite
making brilliant contributions to Maxwell's theory and being awarded 
an honorary doctorate from G\"{o}ttingen University in 1905. The other
``Maxwellians'' didn't make it into the limelight either.

The second big lesson I learned from my encounter with Carl is to
never ignore the Institute for Advanced Study (IAS) in Princeton or the
pround influence of its two first founding permanent members: Hermann
Weyl and Albert 
Einstein. Carl was a student of Salomon Bochner and thrived on all the
mathematics emanating from the IAS. In retrospect, the world seems
quite small. It was Weyl who in 1948 invited Bott to the IAS, it was
Weyl who earlier got de Rham, Kodaira, and Spencer to put Hodge theory
on a rigorous footing, and it was Weyl \cite{W},\cite{Y} and his close
colleague Einstein who were the true curators of the developments
arising from 
Maxwell's theory. It turns out that the ``Whitney forms'' that I was so
fond of have their origins in a 1952 paper of Andr\'{e} Weil called ``Sur
les Theorems de de Rham.'' Clearly, every part of the novel mathematics
I was using could be traced back to the IAS; even if my application of
these ideas to computational electromagnetics was unforeseen. If ever
I was in denial about details, I could check in with Donald Spencer's student,
Robert Hermann, to verify facts\footnote{The ``Advanced Calculus''
 text Nickerson, Spencer, and Steenrod
\cite{NSS} was Princeton-inspired but was never
published. However, it initiated a wave of differential-form
based multi-variable calculus texts in the 1960s. Although it is a
very natural way to bring 
multivariable calculus to its roots in Physics, this wave of texts
didn't catch on. Bott never wrote a text for such an
undergraduate audience and so one can only hypothesize about how he
would have integrated Kirchhoff's laws with Hodge theory and Maxwell's
equations. In retrospect, it took a couple
of decades to get things right and ultimately, the books that
reached out to engineers and physicists most effectively were written by Bott's
close colleagues \cite{BS}, \cite{F}.}. If details were scarce in the
literature, contemplating the influence of Hermann Weyl could help
bring things into focus. 

Enough said about my interactions with Carl Herz. To appreciate Carl Herz'
contributions to harmonic analysis and other fields of mathematics, as
well as the feisty character himself, through the eyes of his
colleagues, the reader is referred to other sources \cite{D},
\cite{H}. Needless to say, when I realized that Carl Herz was behind the 
effort to award Bott an honorary doctorate from McGill, it seemed like
a big piece of the puzzle fell into place. He plays a central role on
almost all correspondence with the university administration on the
matter and in the end, he's the one who presented Bott for the degree
in June of 1987.  

\specialsection*{The convocation presentation of Carl Herz}

\begin{quotation}
Mr. Chancellor,

I have the honour of presenting to you, in order that you may confer
on him the degree of Doctor of Sciences, honoris causa, Professor
Raoul Bott.

Raoul Bott was born in Hungary, but his university education up to the
M. Eng. Was at McGill. He received his B. Eng. from McGill in
1945. After a short stint in the infantry, he continued his studies in
electrical engineering at McGill. The immediate postwar period saw a
great demand for mathematics teachers, and Bott taught calculus here
while studying for his master's degree. In addition he took some
courses from Professor Gilson, then Chair of the Department of
Mathematics. Nevertheless, he remained a student of electrical
engineering until he left McGill to go to Carnegie Tech for his 
doctorate.

Electrical engineering has a close affiliation with what might be
viewed as an abstruse branch of mathematics, algebraic topology,
Professor Bott's specialty. One has only to recall that ``Betti
numbers'', the fundamental numerical invariants of topology, are named
for an Italian electrical engineer, and one can read James Clerk
Maxwell for profound insights into the subject. At an even more
primitive level, circuit theory has always been a source of good
problems for topologists. Bott's earliest work was rather
algebraic. The Bott--Duffin Theorem (1949) on circuit synthesis was
described by a reviewer thus: ``This proof of the realizability of the
driving point impedance without the use of transformers is one of the
most interesting developments in network theory in recent years.'' It
continues to be a much-cited result. This work came shortly after Bott
had obtained a D.Sc. in mathematics.

After the doctorate, Raoul Bott went to the Institute for Advanced
Study in Princeton. He was at the Institute during 1949--1951 and
returned in 1955--1957. He joined the faculty of the University of
Michigan in 1951 where he remained until 1959 when he was invited to
his present academic  home, Harvard, where he is William Caspar
Graustein Professor of Mathematics.

Professor Bott's seminal contributions to mathematics are too
extensive for me to do justice to them here.  Most of his early ideas
seem to have drawn their inspiration from the Calculus of Variations
in its global version known as ``Morse Theory''. Bott applied Morse
Theory in an unexpected and striking way. Over a long period he,
together with his various collaborators, worked out the topology of
Lie groups and symmetric spaces. One must mention the Bott Periodicity
Theorem which brought some order to the chaos of homotopy theory. He
went on to study fixed point theorems and their application to other
branches of mathematics including differential equations.  Most
recently Bott has been working on applications of topology and
geometry to the Yang-Mills equations in quantum field theory.

For his achievements, Bott was awarded the Veblen Prize of the
American Mathematical Society in 1964.

In addition to his purely scientific accomplishment, Raoul Bott
stimulates all those who are about him. He is one of the best and most
exciting expositors of mathematics I have had the privilege to listen 
to.

Mr, Chancelor, McGill can take great pride in honoring this year, as
it did last\footnote{Raoul Bott, Jim Lambek and Louis Nirenberg all
  graduated from McGill in 1945 and Nirenberg was awarded an honorary
  doctorate from McGill in 1986.}, another of its graduates who stand
in the forefront of mathematics of the twentieth century. 
\par
The Eleventh Day of June, Nineteen Hundred and Eighty-seven

\par
Carl Herz

\par
Professor of Mathematics and Statistics.
\end{quotation}

Given my encounters with Prof. Herz, his correspondence with the McGill
administration, and his encyclopedic breadth, it is clear that he
played a central
role in the case for the honorary degree. The masterful presentation
of Carl Herz shows how a broad perspective can lead to a
reorganization of knowledge that lets the likes of Paul Dirac and Eugene
Wigner move from Engineering to Physics, and the likes of Raoul Bott,
Solomon Lefschetz, John
Milnor, and Donald Spencer move from Engineering to Mathematics. On
the other hand, since Prof. Herz always enjoyed an argument (in the
very best sense of the word!), I'll take the liberty to make a
qualification and perhaps an eleboration. 

The qualification I might add is that Enrico
Betti was clearly not an electrical engineer but a mathematician.
Indeed, Betti made contributions to both Elasticity theory and
Electromagnetism, and Maxwell does indeed cite Betti's work in his
treatise, but he was a mathematician. Betti, like the entire school of
Italian Algebraic Geometry, was highly influenced by Riemann and
topological ideas. However, the level of rigor in 19th century Italy was lax by
modern standards and so his influence on current mathematical research
may seem far removed. I revere Carl's
respect for historical detail, and I'll refrain from calling his
labeling Betti as an Electrical Engineer as a mistake. Rather I'd say
Betti, like Gauss, Riemann, and Vito Volterra, had broad
interests, and that Prof. Herz suppressed the pedant
in himself and took some license in his interpretation of
history. 

\section*{The hidden hand of Hermann Weyl}

The role of Hermann Weyl in getting mathematics off
the ground in the earliest days of the IAS is now well documented \cite{B}.
What I find fascinating is the first and fateful encounter between
Hermann Weyl  
and Raoul Bott. The encounter has a lot to do with interplay between electrical
circuit theory, the early days of algebraic topology, and the
perception of topology. The presentation of Carl Herz leaves out
a lot of detail, much as a movie based on a book has to
forgo a lot of detail. Given the encyclopedic knowledge of Carl Herz,
it is tempting to speculate on what he could have put into a longer
presentation.

Bott told the story of his first encounter with Hermann Weyl many
times, emphasizing different aspects and different amounts of
detail. See for example\cite{Bo88}. I like
the following rendition of the basic facts. During his grad student
days as a student of Richard Duffin at Carnegie Tech, Bott played a
large role in organizing the department colloquium. Being fluent in
German, Hungarian, and Slovakian he would have an edge over other grad
students in terms of ``chatting up'' foreign-born visitors. When
Hermann Weyl visited, they were introduced, and Bott immediately began
to tell Weyl of his thesis
work\cite{BDu}. There are interesting aspects of his thesis which
predate both the Bott--Duffin Synthesis procedure and Wang
algebras\cite{D59}. One key aspect is the 
``impedance potential'' and how it defines a generalized
inverse of a matrix. Of course the Moore--Penrose axioms for a generalized
inverse were only formulated in the 1950s and so Bott does not use the
term.\footnote{See Chapter 2 and Appendix A of
 Ben-Israel and Greville's book \cite{BIG} for an exposition that
 puts the Bott-Duffin constrained inverse in the context of
 generalized inverses, and for putting the ``Moore'' of
``Moore-Penrose'' in historical perspective.} Early on, he and Duffin
called it a ``constrained inverse''\cite{BDu}, and in expository talks
Bott later described it in terms of orthogonal projections in a complex
(i.e. Hodge theory).  It turns out that his
impedance potential is a determinant which is intimately related to
what graph theorists call a ``Matrix-tree formula''-- a result that
goes back to Kirchhoff and was used in Maxwell's treatise. The logarithmic
derivative of the impedance potential with respect to
branch impedances gives a generalized inverse. When Bott
explained the formalism and associated results to
Hermann Weyl, Weyl grasped that Bott-Duffin synthesis was indeed a
contribution to network synthesis, but that the connection between
Hodge theory and Kirchhoff's laws was not. He pointed Bott to some
papers connecting Kirchhoff's laws to topology which he wrote in the
early 1920s\footnote{More precisely, Weyl's papers dealt with
  Kirchhoff's laws \cite{W23a} and combinatorial topology
  \cite{W23b}.}. Needless to say, Bott was invited to the IAS, but Bott 
felt a bit deflated about the Hodge theoretic aspect and that Weyl saw
it concretely in Kirchhoff's work. 

To be fair to Bott, we have to ask why these papers of Hermann Weyl
were so obscure. Gottingen was very closely tied to the technological
aspects of Maxwell's theory, and so why were these two papers as
obscure as Maxwell's periphractic numbers? What was the point Weyl was
trying to make? To give some insight, a digression is in order.

In the Winter of 2005 I spent a month in the math department at the
ETH in Zurich while on sabbatical. When I arrived, my host gave me a 
choice of offices: a huge office with a stunning view of Zurich belonging
to a colleague on Sabbatical, or a very small empty office in the back of
the building where ``pure mathematicians hide their guests.'' I told
my host that I wanted the freedom to ``spread out,'' and that I felt
more comfortable in the small back office. He was perplexed but
obliged. It turns out that my cozy office was next to that of Beno
Eckmann. On the centennary of Einstein's golden year, Zurich celebrated
Eckmann as the last person in the city who had 
personal contact with Einstein! Since Hermann Weyl was the the head of
the ETH mathematics department in the 1920s, I naturally wanted to
pick Beno's brain for anecdotes. Not wanting to mess with his work
habits, I planned to chat him up while he was a sitting target.  

In the hallway outside our offices was a high-tech espresso
machine and every morning Beno would take a break to sit and enjoy an
espresso outside our offices. The first day, I ``coincidentally''
joined him and he related wonderful anecdotes from 1950-1955, after
Hermann Weyl retired from the IAS study, resettled in Zurich, and
frequented the department. (I was sufficently impressed that when I
returned to Boston, I contacted the editors of the Notices of the
AMS and a year later the anecdotes appeared in print\cite{EWZ}).  The
next day I resolved
to ask Beno a question which I didn't think any living person could
answer. Little did I know that he had written a paper on the subject \cite{E}
and had a definite opinion on every nuance I could ask him to
elaborate on! 

The conversation went something like this:

RK: Beno, there is something I really don't understand about Hermann
Weyl.

BE: What is it?

RK: Well, in his collected works, there are are two papers about
electrical circuit theory and topology dating from 1922/3. 
They are written in Spanish and published in an
obscure Mexican mathematics journal. They are also the only papers he ever
wrote in Spanish, the only papers published in a relatively obscure
place, and just about the only expository papers he ever wrote on
algebraic topology. It would seem that he didn't want his
colleagues to read these papers.

BE: Exactly!

RK: What do you mean?

BE: Because topology was not respectable!

RK: Why was toplogy not respectable?

BE: Hilbert!

RK: Hilbert?

BE: Just look at his 23 problems from 1900. Do you see anything to do
with combinatorial group theory or topology? No!

RK: Why?

BE: Poincar\'{e}!\footnote{Sarkaria\cite{S} has given a modern executive
 summary of Poincar\'{e}'s work in Topology.}

RK: What did Hilbert think of Poincar\'{e}'s work on toplogy?

BE: Poincar\'{e} would write a huge paper on Analysis Situs. Half of it
would be completely wrong! So, he'd write another huge paper trying to
correct the first, but it would be half wrong! And so he'd write a
third paper, but it would be half wrong. And so on... deuxieme
complement, troisieme, quatrieme, cinquieme,.. and in the end what
did we get? Dubious results and conjectures! Hilbert didn't think this
was mathematics!

RK: So why did Hermann Weyl write these papers?

BE: He wanted to take stock of the honest results and
reorganize them using a more modern abstract algebraic approach. Emmy
Noether and others were doing interesting things in algebra and he had a need to
write these papers for himself. These papers also contain some new
results like the signature of a 4-d manifold.

Beno went on to portray Hilbert as a bit of a reactionary
figure, around which Hermann Weyl had to tip-toe. However, if Weyl
wanted an opportunity to move things forward, it came in 1930 when Weyl
succeeded Hilbert upon his retirement from G\"{o}ttingen. Although he was only
at the helm from 1930 until he fled the Nazis in 1933, during this
very brief time German topology flowered in the hands of Emil Artin, Kurt
Reidemeister, and others. According to Beno Eckmann, Weyl made a
historic decision in 1930 which was highly controversial at the time,
but ultimately 
vindicated: he appointed Heinz Hopf, a young researcher and relatively
unseasoned, as his successor at ETH. 

If Beno's historical perspective is taken superficially, there is a
temptation to
suspect there was some lasting disagreement between Weyl and
Hilbert. However, one only needs to read Weyl's masterful summary of
Hilbert's work \cite {W44} to realize that both men held themselves to
the highest standards. In a sense, every time Hilbert or Poincar\'{e}
dug their heels in, Weyl found an opportunity to move mathematics
forward. Algebraic topology may 
be one example and the continuum hypothesis may be another; perhaps
the best example is the fact that Kurt G\"{o}del was among the first four
hires at the IAS, but unemployable in Europe.

What did Carl know? We can only speculate. I can only say that he is
one of the many people who impressed upon me the importance of having
a historical perspective when reconciling algebraic topology with its
applications.

\section*{From History to Bott's reconciliation with it}

The historical details discussed so far predate 1952. The next eight
years would usher in the revolution in homotopy theory brought on by
Serre's thesis, CW complexes which tie Morse theory to homotopy
theory, Bott periodicity, generalized cohomology theory such as
K-theory and Rene Thom's cobordism theory, and the reformulation of
generalized cohomology theories
in terms of spectra. Beno Eckmann pointed out to me that in the 1950s,
those in Zurich who dismissed Hermann Weyl as an old man were
favorably stunned by his
summary of the work of Kunihiko Kodaira and Jean-Pierre Serre on the
occasion of their being awarded the Fields Medal in 1954 \cite {W54}. 
Nonetheless, contrasting Weyl's presentation of the work of the two Fields
medalists, it is apparent that he was challenged by the homotopy
theoretic world Bott had entered into, even if he did a lot to unleash
the homotopy theoretic perspective. Enough said; it is time to leave
threads of mathematical history
and experience another view of history: 

\specialsection*{The convocation address of Raoul Bott}
\begin{quotation}
Mr. Chancellor, Mr. Chairman of the Board-- my dear fellow graduates:

Congratulations to you- class of '87! You look splendid! I think you
wash more behind the ears than your American cousins at Harvard do.

It is nice to get a degree, isn't it? Of course you only had to work
hard for four years or so to get yours, while it took me over forty
years to get mine. And presumably you have paid for yours, while I am
paying for mine at this moment by being here on this platform, making
a fool of myself. 

But, there is really nothing like one's first degree.  And what I
loved especially about my Bachelor of Engineering was that an iron
ring (from a fallen bridge) came with it. I hope this tradition
continues, so that at least you engineers, can contrive -- as I did --
to display it on every occasion. It is a marvelous way of starting a
conversation and at the same time lets one know that you have
``graduated.'' So my first admonition to you is: ``Flaunt your degree in
front of the whole world!''  For a few weeks enjoy it to the hilt! The
real world will rein you in soon enough.

Of course the people who enjoy your degree most are your parents. So
by all means--- here comes my second admonition --- Get yourselves some
children, in time for degree-harvesting when you are still in your
forties! (That way you might also have time to repay the loans before
you die.)

But let me tell you now a little bit about the good old days, just to
keep some sort of historical perspective in a society, whose customs
change at such a rate that the last forty years most probably
represent two hundred uninflated ones.

First of all I must tell you that, beautiful as your Campus is today,
it used to be even more so in 1941. There were lawns to stretch out
on, there was even a tennis court by the Redpath library! There was so
much space and such a fine line of proportion! And there were no
skyscrapers! (On the other hand, the area around McGill was very 
rundown. And I see that our ``greasy  spoon'' has now flowered into a
pizza joint.)

Classes were small and some of my professors wore robes, as we are
now, to teach in. They billowed and flowed delightfully with each
step.  These gowns were usually torn and completely covered with dust;
still they added to the performance. I remember that later when I
had my own calculus class to teach -- the veterans had returned in huge
numbers in the fall of 1945 and the Math. Department had pressed a
lowly engineer into service to meet the demand -- my dear friend and
mentor, Prof. McLennan -- the Socrates of our campus -- lent me his
well-weathered gown. ``Try it in your class'', he said, with a twinkle
in his eye. Well, the class of course guffawed at first, but then
actually settled down to work in a more businesslike manner than usual.

Possibly it was this ballet-like aspect of the lectures that kept me
going to classes very diligently in the beginning. However, this epoch
of my life came to an end in short order after one of my roommates in
our boarding house on Durocher called me in for a serious talk. Elwood
Henneman was his name and he continues as a dear friend and colleague
at Harvard. Elwood, with the full authority of a first year medical
student and a Harvard B.A. warned me of the danger of being addicted
to classes. ``Never become a slave to them,'' he declared; ``do the bulk
of your thinking on your own!''  This point of view made immediate
sense; and thereafter it was safest to look for me in the Music
room. Usually in the company of my dear friend, Walter Odze, and much
later on also with my wife to be.

How much this had to do with my falling grades I don't know, but in any
case I did manage to always to ``beat the Dean'' as we used to say. Do
you know what I mean? (But it sounds like good fun in any case,
doesn't it?) Well, one of the endearing procedures of our Alma Mater
at that time was that they didn't divulge our final grades until
August -- I think. Then suddenly your name was printed in the Gazette --
if you passed, that is -- and with an asterisk if you flunked one
course, etc. On the same day the names of all the passing students
were listed on a billboard in linear order of merit. Those who had
failed were not on the list, and at the bottom of this terrifying
document came the signature of the dean! Hence the expression of
beating the Dean if one got through.

But speaking of Deans and advice, let me tell you about one McGill
Dean who in his own inimitable way gave me the best advice of my
life. These were the war years and in '45 right after graduation, I
joined up in the Canadian infantry and was being trained for combat in
Japan. After three months in basic training the atomic bomb was
dropped on Hiroshima and Nagasaki, the war ended abruptly and my
fellow recruits and I thereby suddenly and miraculously reprieved -- in
this unbelievable and terrible manner.

Of course, the one great advantage of being in the army is that one
has no career problems whatsoever! Hence the doubts I had about my
vocation in engineering were completely submerged by my efforts to
keep out of the Sergeant's hair.  But when in October I found myself
back in the Engineering Department, where they had very kindly let me
return for a Master's Degree on -- as you can imagine -- very short
notice, the old doubts flared up again and I was in a quandary about
what to do.

It was sometime in `46 then that I presented myself at Dean Thompson's
office and asked him whether he could see his way to putting me
through medical school.  (On the Jewish side of my family they always
did say: ``chutzpah he does not lack''.) And Dean Thompson was quite
encouraging at first. ``We need scientifically trained doctors'', he
said. ``But'', he continued, ``first tell me a little about yourself''. It
was at this point that our interview started to go sour. No, I never
enjoyed Biology much. No I hated dissecting frogs. Alas, Botany bored
me and I had little use for Chemistry! After this sorry litany,
Dr. Thompson surveyed me and the situation for a while, pipe in hand,
and lost in thought. ``Is it maybe that you want to do good for
humanity'', he said at last.  I hemmed and hawed in my seat, but before
I had time to say anything he came out with: ``Because they make the
lousiest doctors!''

Well, that was it for me -- and you must admit that it explains a lot
of things, doesn't it? In any case, I got up and as I went to the
door, I thought to myself: well you (explicative deleted) if that is
how the land lies, then I will simply do what I like best: ``I will
become a Mathematician. Put that in your pipe and smoke it!''

I tell you all this only in part as a jest.  I would also like it to
be a word of encouragement to those of you who, degree in hand, still
are not quite certain of your path. May you also be blessed with a
counselor with such diagnostic skills and such a knack for putting you
on the right course.

But my time is up! Still I cannot resist a serious word.  Over my
McGill days the War hung like an everpresent black cloud, subtly
affecting every aspect of our lives. For you in the nuclear age the
cloud is, thank God, farther away, but potentially much, much
darker. These things you will have to live with and somehow hope to
conquer. But for this road I know of no better advice than my friend
Elwood's - ``Do your own thinking''.

In our more immediate lives we are also beset today more than ever
before, with show, with image, with jargon; and here again -- to pick
one's way through this quagmire, there is no better exhortation than:
Be your own man; be your own woman. For then, I am confident, you will
never confuse fashion with substance, heroes with the people who
depict them on the tube, computers with people, Science with virtue,
or wealth with happiness.

Yes, may God bless you and may you be joyfully and productively
yourselves - but may you also be ultimately servants of a larger and
an all-encompassing benign world  view. And that is really no more
than what I take to be the correct reading of Dr. Thompson's advice to
me forty years ago. Only remember that the concerns of your generation
must even transcend those off our human family. They must embrace
every aspect of life itself on this deeply troubled, but magnificent
and magical planet of ours.
\par
Raoul Bott
\end{quotation}

\section*{On Iron Rings and other aspects of the Convocation Address}

Bott's convocation address sets the stage for my fascination with
his struggle to reconcile old and new world sensibilities.  He was
clearly the same age as my parents and 
like my parents, the disruption of his adolescence by Hitler, Stalin,
and the events in the Europe of his youth was a traumatic
experience and a profound education even if they didn't view it that
way at the time\cite{BoD}. His refugee experience was a stark
contrast to the way kids grew up in North America during the decades
after WWII. This was evident in his sense of humor and in the way he
handled those who did not choose their words correctly. The
testimonies of his students in this volume attest to this.  The war
years and the economic turmoil that preceded it left him with little
tolerance for the dangerous comforts of self-preservation. The
Engineer's iron ring along with its uniquely Canadian
origin seem to frame some of his advice on responsibility and independent
thinking. Reading the convocation speech, and recalling our
interaction after the first class I audited at Harvard, it is apparent
that Bott 
had a much more profound appreciation of, and respect for, the iron
ring than did students he was addressing. One cannot do
justice to the topic here\footnote{These days one has to
  look for stainless steel if one wants to spot an ``Iron
  Ring''---the iron of the early rings used to be eaten away by
  sweat and was soon replaced.}, but it is useful to connect a few key ideas
to events and sensibilities of other times. 

Engineering is full of trade-offs, and the story of the iron ring is
about the interface between technical trade-offs, ethics and ambition. One
engineering trade-off is between the theoretical effort that goes into
designing something without making a physical model, and the
willingness to build prototypes and make mistakes. If one were to
design a paper clip, one would make many prototypes in order to see
``what works.'' On the other hand, if one were building a bridge, one
would like to avoid disasters, and the development of a theoretical
model with predictive properties is in order. In the case of bridge
building, especially new designs, one is very cautious because public
confidence is paramount. However, there have been many bridge
disasters and contrary to what one might naively expect, they usually
do not involve new designs! They usually involve refined designs which
take into account that earlier designs were overly cautious, too
costly, and less than ambitious. The temptations involved are quite
universal and are not restricted to bridges; one can send
up space shuttles routinely without an accident, but when one decides that
the rules for launch are overly cautious in light of an
opportunity to make some ``State of the Union Address''
spectacular, strange things can happen --- just like when the detailed
properties of O-rings were purposely ignored in the lead-up to the Challenger
disaster. Similarly, the dismissing of foam
impacts during launch as routine in the lead-up to the Columbia
Shuttle disaster underlines the vigilance reqired to make
complicated things work. These days we have ``financial engineering'' and
computer models 
so predictive that there is a temptation to lose track of underlying
assumptions and to consider the regulation of investors as
unimaginative and
cumbersome- here again, strange things can happen when ambition trumps
regulation. Engineering disasters create
teachable moments and they are very well documented when the stakes are
high. In the case of bridges, the spectacular disasters have been
studied and categorized, and scholarly books such as the one
by Petroski\cite {P} have been written. 

 Chapter three of Henry Petroski's book details the Quebec City bridge
 disaster(s) that lead to the iron ring worn by Canadian Engineers. It has a
 lot of detail on the New York based construction firm, the details of
 the bridge, the ignoring of warning signs and the 75 people
 killed in the first disaster. A key ingredient in this August 1907
 disaster was an attempt to redesign the bridge during construction in
 order to ensure it broke a world record. There was a second disaster
 in 1916 during the construction of the  
 redesigned bridge which killed 13 workers. In all, a total of 89
 workers were killed in the construction of the bridge. The completion
 of the 1800 foot span of the Quebec City bridge in 1917 made it 
 the largest cantilever brige in the world and vindicated the concept of
 the cantilever bridge for a mix of rail and automobile traffic. However,
 worldwide, no other major cantilever 
 bridge was completed until the 1930s. To this day the Quebec City bridge has 
 the longest span of any cantilever bridge -- other bridge
 designs are used for longer spans. The original iron rings were made
 of the collapsed bridge's iron as a reminder of the stupidities
 engineers are capable of, and as a reminder of the engineer's
 responsibility to society -- soon after the rings were made of
 stainless steel.  

Rudyard Kipling was the recipient of the 1907 Nobel Prize in
literature and lived in Bratleboro Vermont for a few years in the
1890s. It was in Vermont, between Quebec City and the home of the bridge's
architect in New York, that he wrote ``The Jungle Book.'' It was also in
1907, following the first Quebec City bridge disaster, that he wrote a
poem called ``The Sons of Martha''.\footnote{See, for instance:
  "Kipling: A Selection of His Stories and Poems" by 
John Beecroft, (in two volumes), Doubleday 1956. In Vol. II, The Sons
of Martha appears on p 451.} It is inspired by the Gospel of Luke
(10:38-42) and forms the basis of the original iron ring
ceremony which Kipling was comissioned to write. The original Canadian
ceremony was called ``The Ritual of the Calling of an Engineer'' and
it was first performed in 1922. These days the original ceremony with
its Biblical references is considered noninclusive. The iron ring 
ceremony, first performed in the United States in 1970, is centered
around ``The obligation of the engineer'' which is devoid of Biblical
references.\footnote{Kipling lived after ``the days of wooden ships and
 iron  men'', and in the peak  of the English Empire. Of the 89 workers killed 
in the two Quebec City bridge disasters, it appears that 33 were
Mohawk steel workers form the Kahnawake reserve just outside of
Montreal, creating 24 widows and numerous fatherless children. The
Mohawk workers were well adapted to heights and were the ``high tech''
workers of their time. I have yet to see this aspect arise in the
context of a modern iron ring cremony, or in the world view of
Kipling's time.} In Bott's commencement
address he is clearly referring to the the original ceremony and I
recommend a read through Kipling's poem to appreciate the iron ring as
Bott would have been introduced to it.   

This concludes my musings on a refugee as an engineering
student in Montreal, his metamorposis into a topologist, and a new
world success story who was ultimately awarded an honorary
doctorate by his alma mater.

\bibliographystyle{amsalpha}

\end{document}